\documentclass{amsart}
\usepackage{amsmath,amssymb}
\usepackage{yhmath}

\newcommand{\bdis}{\begin{displaymath}}
\newcommand{\edis}{\end{displaymath}}
\newcommand{\be}{\begin{equation}}
\newcommand{\ee}{\end{equation}}
\newcommand{\mbb}{\mathbb}
\newcommand{\mcal}{\mathcal}

\newcommand{\vp}{\varphi}

\newcommand{\zf}{\zeta\left(\frac{1}{2}+it\right)}
 
\newcommand{\zfn}{\zeta\left(\frac{1}{2}+it_\nu\right)}

\newcommand{\FR}{\frac{x^n+y^n}{z^n}}

\theoremstyle{definition}

\theoremstyle{remark}
\newtheorem{remark}[]{Remark}

\newtheorem*{mydef11}{{\bf Theorem 1}}

\newtheorem*{mydef12}{{\bf Theorem 2}}

\newtheorem*{mydef13}{{\bf Theorem 3}}

\newtheorem*{mydef14}{{\bf Theorem 4}}

\newtheorem*{mydef41}{{\bf Corollary 1}}

\newtheorem*{mydef51}{{\bf Lemma 1}}

\newtheorem*{mydef52}{{\bf Lemma 2}}

\newtheorem*{mydef81}{{\bf Property 1}}

\newtheorem*{mydef82}{{\bf Property 2}}

\newtheorem*{mydefII1}{\bf Consequence 1}

\numberwithin{equation}{section}

\begin{document}

\title[Jacob's ladders, new equivalent \dots]{Jacob's ladders, new equivalent of the Fermat-Wiles theorem generated by certain cross-breed of Ingham and Heath-Brown formula (1979) and some chains of equivalents}

\author{Jan Moser}

\address{Department of Mathematical Analysis and Numerical Mathematics, Comenius University, Mlynska Dolina M105, 842 48 Bratislava, SLOVAKIA}

\email{jan.mozer@fmph.uniba.sk}


\begin{abstract}
In this paper we obtain a new equivalent of the Fermat-Wiles theorem based on a kind of cross-bred of Ingham and D. R. Heath-Brown formula. Further, we prove the existence of infinite set of finite chains of a kind of equivalent expressions of mathematical analysis. 
\end{abstract}
\maketitle

\section{Introduction} 

\subsection{} 

In this paper, that continues our series \cite{8} -- \cite{16}, we use our method of constructing new functionals on:  
\begin{itemize}
	\item[(A)] the Ingham's and D.R. Heath-Brown's formula  (see \cite{2}, \cite{3}, p. 129)
	\be \label{1.1} 
	\begin{split}
	& \int_0^T\left|\zf\right|^4{\rm d}t=T\sum_{s=0}^4 a_s\ln^{4-s}T+\mcal{O}(T^{7/8+\epsilon}), \\ 
	& a_0=\frac{1}{2\pi^2},\ T\to\infty, 
	\end{split}
	\ee 
	\item[(B)] and our quotient-formula (see \cite{16}, (3.11))
	\be \label{1.2} 
	\frac{\int_{\overset{r-1}{T}}^{\overset{r}{T}}|\zf|^2{\rm d}t}{\int_{\overset{r-1}{T}}^{\overset{r}{T}}|\zeta(\sigma+it)|^2{\rm d}t}=\frac{1}{\zeta(2\sigma)}\ln\overset{r}{T}+\mcal{O}(1) 
	\ee 
	for every fixed small $\epsilon>0$ and every fixed 
	\bdis 
	\sigma\geq\frac 12+\epsilon. 
	\edis 
\end{itemize} 

\begin{remark}
Let us remind, that the original formula 
\be \label{1.3} 
\int_0^T\left|\zf\right|^4{\rm d}t=\frac{1}{2\pi^2}\ln^4T+\mcal{O}(T\ln^3T) 
\ee 
with the valid first term belongs to Ingham (1926). 
\end{remark} 

\subsection{} 

We obtain in this paper a new functional as a kind of cross-breed of the formulas (\ref{1.1}) and (\ref{1.2}), namely 
\be \label{1.4} 
\begin{split} 
& \lim_{\tau\to\infty}\frac{1}{\tau}
\left\{
\int_{\frac{2\pi^2}{\zeta^4(2\sigma)}x\tau}^{[\frac{2\pi^2}{\zeta^4(2\sigma)}x\tau]^1}|\zeta(\sigma+it)|^4{\rm d}t
\right\}^4\times 
\left\{
\int_0^{\frac{2\pi^2}{\zeta^4(2\sigma)}x\tau}\left|\zf\right|^4{\rm d}t
\right\}\times \\ 
& \left(\sum_{s=0}^4 c_s
\left\{ 
\int_{\frac{2\pi^2}{\zeta^4(2\sigma)}x\tau}^{[\frac{2\pi^2}{\zeta^4(2\sigma)}x\tau]^1}\left|\zf\right|^2{\rm d}t
\right\}^{4-s}\right)^{-1}\times \\ 
& \left( \left\{
\int_{\frac{2\pi^2}{\zeta^4(2\sigma)}x\tau}^{[\frac{2\pi^2}{\zeta^4(2\sigma)}x\tau]^1}|\zeta(\sigma+it)|^2{\rm d}t
\right\}^s \right)^{-1}=x, 
\end{split} 
\ee 
for every fixed $x>0$ and $\sigma\geq\frac 12+\epsilon$, where 
\bdis 
c_0=1,\ c_s=2\pi^2\{\zeta(2\sigma)\}^{-s}a_s,\ s=1,2,3,4. 
\edis 

\begin{remark}
The new equivalent of the Fermat-Wiles theorem follows from (\ref{1.4}) in the special case of Fermat's rational  
\bdis 
x\to\FR,\ x,y,z,n\in\mbb{N},\ n\geq 3. 
\edis 
\end{remark}

\subsection{} 

Further, we also give in this paper an infinite set of finite chains of asymptotic equivalent expressions (integrals, compositions of integrals, sums and functions) corresponding to our first functional (see \cite{9}, (4.6))
\be \label{1.5} 
\lim_{\tau\to\infty}\frac{1}{\tau}\int_{\frac{x}{1-c}\tau}^{[\frac{x}{1-c}\tau]^1}\left|\zf\right|^2{\rm d}t=x, \ x>0.  
\ee  
For example, a chain corresponding to the value $x=1-c$ is the following one: 
\be \label{1.6} 
\begin{split} 
& \int_\tau^{\overset{1}{\tau}}\left|\zf\right|^2{\rm d}t\sim 
\left\{
\int_{\frac{2\pi^2(1-c)}{\zeta^4(2\sigma)}\tau}^{[\frac{2\pi^2(1-c)}{\zeta^4(2\sigma)}\tau]^1}|\zeta(\sigma+it)|^2{\rm d}t
\right\}^4\times \\ 
& \int_0^{\frac{2\pi^2(1-c)}{\zeta^4(2\sigma)}\tau}\left|\zf\right|^4{\rm d}t \times \\ 
& \left(
\sum_{s=0}^4c_s\left\{\int_{\frac{2\pi^2(1-c)}{\zeta^4(2\sigma)}\tau}^{[\frac{2\pi^2(1-c)}{\zeta^4(2\sigma)}\tau]^1}\left|\zf\right|^2{\rm d}t\right\}^{4-s}\times \right. \\ 
& \left. \left\{\int_{\frac{2\pi^2(1-c)}{\zeta^4(2\sigma)}\tau}^{[\frac{2\pi^2(1-c)}{\zeta^4(2\sigma)}\tau]^1}|\zeta(\sigma+it)|^2{\rm d}t\right\}^s
\right)^{-1} \sim \sum_{\tau< n\leq \overset{1}{\tau}}d(n)\sim\sum_{\tau<t_\nu\leq\overset{1}{\tau}}\pi\zfn \\ 
& \sim \ln\left\{\frac{\Gamma(\overset{1}{\tau})}{\Gamma(\tau)}\right\}\sim \int_1^{\frac{1-c}{\zeta(2\sigma)}\tau}|\zeta(\sigma+it)|^2{\rm d}t\sim 
\int_0^{\frac{1-c}{\bar{c}_l}\tau}|S_1(t)|^{2l}{\rm d}t,\ \tau\to\infty. 
\end{split} 
\ee

\section{Jacob's ladders: notions and basic geometrical properties}  

\subsection{} 

In this paper we use the following notions of our works \cite{4} -- \cite{7}: 
\begin{itemize}
\item[{\tt (a)}] Jacob's ladder $\vp_1(T)$, 
\item[{\tt (b)}] direct iterations of Jacob's ladders 
\bdis 
\begin{split}
	& \vp_1^0(t)=t,\ \vp_1^1(t)=\vp_1(t),\ \vp_1^2(t)=\vp_1(\vp_1(t)),\dots , \\ 
	& \vp_1^k(t)=\vp_1(\vp_1^{k-1}(t))
\end{split}
\edis 
for every fixed natural number $k$, 
\item[{\tt (c)}] reverse iterations of Jacob's ladders 
\be \label{2.1}  
\begin{split}
	& \vp_1^{-1}(T)=\overset{1}{T},\ \vp_1^{-2}(T)=\vp_1^{-1}(\overset{1}{T})=\overset{2}{T},\dots, \\ 
	& \vp_1^{-r}(T)=\vp_1^{-1}(\overset{r-1}{T})=\overset{r}{T},\ r=1,\dots,k, 
\end{split} 
\ee   
where, for example, 
\be \label{2.2} 
\vp_1(\overset{r}{T})=\overset{r-1}{T}
\ee  
for every fixed $k\in\mbb{N}$ and every sufficiently big $T>0$. We also use the properties of the reverse iterations listed below.  
\be \label{2.3}
\overset{r}{T}-\overset{r-1}{T}\sim(1-c)\pi(\overset{r}{T});\ \pi(\overset{r}{T})\sim\frac{\overset{r}{T}}{\ln \overset{r}{T}},\ r=1,\dots,k,\ T\to\infty,  
\ee 
\be \label{2.4} 
\overset{0}{T}=T<\overset{1}{T}(T)<\overset{2}{T}(T)<\dots<\overset{k}{T}(T), 
\ee 
and 
\be \label{2.5} 
T\sim \overset{1}{T}\sim \overset{2}{T}\sim \dots\sim \overset{k}{T},\ T\to\infty.   
\ee  
\end{itemize} 

\begin{remark}
	The asymptotic behaviour of the points 
	\bdis 
	\{T,\overset{1}{T},\dots,\overset{k}{T}\}
	\edis  
	is as follows: at $T\to\infty$ these points recede unboundedly each from other and all together are receding to infinity. Hence, the set of these points behaves at $T\to\infty$ as one-dimensional Friedmann-Hubble expanding Universe. 
\end{remark}  

\subsection{} 

Let us remind that we have proved\footnote{See \cite{8}, (3.4).} the existence of almost linear increments 
\be \label{2.6} 
\begin{split}
& \int_{\overset{r-1}{T}}^{\overset{r}{T}}\left|\zf\right|^2{\rm d}t\sim (1-c)\overset{r-1}{T}, \\ 
& r=1,\dots,k,\ T\to\infty,\ \overset{r}{T}=\overset{r}{T}(T)=\vp_1^{-r}(T)
\end{split} 
\ee 
for the Hardy-Littlewood integral (1918) 
\be \label{2.7} 
J(T)=\int_0^T\left|\zf\right|^2{\rm d}t. 
\ee  

For completeness, we give here some basic geometrical properties related to Jacob's ladders. These are generated by the sequence 
\be \label{2.8} 
T\to \left\{\overset{r}{T}(T)\right\}_{r=1}^k
\ee 
of reverse iterations of of the Jacob's ladders for every sufficiently big $T>0$ and every fixed $k\in\mbb{N}$. 

\begin{mydef81}
The sequence (\ref{2.8}) defines a partition of the segment $[T,\overset{k}{T}]$ as follows 
\be \label{2.9} 
|[T,\overset{k}{T}]|=\sum_{r=1}^k|[\overset{r-1}{T},\overset{r}{T}]|
\ee 
on the asymptotically equidistant parts 
\be \label{2.10} 
\begin{split}
& \overset{r}{T}-\overset{r-1}{T}\sim \overset{r+1}{T}-\overset{r}{T}, \\ 
& r=1,\dots,k-1,\ T\to\infty. 
\end{split}
\ee 
\end{mydef81} 

\begin{mydef82}
Simultaneously with the Property 1, the sequence (\ref{2.8}) defines the partition of the integral 
\be \label{2.11} 
\int_T^{\overset{k}{T}}\left|\zf\right|^2{\rm d}t
\ee 
into the parts 
\be \label{2.12} 
\int_T^{\overset{k}{T}}\left|\zf\right|^2{\rm d}t=\sum_{r=1}^k\int_{\overset{r-1}{T}}^{\overset{r}{T}}\left|\zf\right|^2{\rm d}t, 
\ee 
that are asymptotically equal 
\be \label{2.13} 
\int_{\overset{r-1}{T}}^{\overset{r}{T}}\left|\zf\right|^2{\rm d}t\sim \int_{\overset{r}{T}}^{\overset{r+1}{T}}\left|\zf\right|^2{\rm d}t,\ T\to\infty. 
\ee 
\end{mydef82} 

It is clear, that (\ref{2.10}) follows from (\ref{2.3}) and (\ref{2.5}) since 
\be \label{2.14} 
\overset{r}{T}-\overset{r-1}{T}\sim (1-c)\frac{\overset{r}{T}}{\ln \overset{r}{T}}\sim (1-c)\frac{T}{\ln T},\ r=1,\dots,k, 
\ee  
while our eq. (\ref{2.13}) follows from (\ref{2.6}) and (\ref{2.5}). 

\section{Transformation of the formula (\ref{1.1})} 

\subsection{} 

First, we rewrite the sum on the right-hand side of (\ref{1.1}) into the following form 
\be \label{3.1} 
\begin{split}
& =\frac{1}{2\pi^2}\sum_{s=0}^4b_s\ln^{4-s}T, \\ 
& b_0=1,\ b_s=2\pi^2a_s,\ s=1,2,3,4, 
\end{split}
\ee 
and next, we use the substitution\footnote{See (\ref{1.2}) with $r=1$.}
\be \label{3.2} 
\begin{split}
& \ln T=\zeta(2\sigma)\frac{\int_T^{\overset{1}{T}}|\zf|^2{\rm d}t}{\int_T^{\overset{1}{T}}|\zeta(\sigma+it)|^2{\rm d}t}\left\{1+\mcal{O}\left(\frac{1}{\ln T}\right)\right\}= \\ 
& \zeta(2\sigma)\times C \times \left\{1+\mcal{O}\left(\frac{1}{\ln T}\right)\right\},\ T\to \infty, 
\end{split}
\ee 
with the shortened notation $C$ for the quotient, in our formula (\ref{3.1}) and obtain 
\be \label{3.3} 
\begin{split}
& =\frac{1}{2\pi^2}\sum_{s=0}^4b_s\{\zeta(2\sigma)\}^{4-s}C^{4-s}\left\{1+\mcal{O}\left(\frac{1}{\ln T}\right)\right\}\\ 
& = \frac{\zeta^4(2\sigma)}{2\pi^2}\sum_{s=0}^4b_s\{\zeta(2\sigma)\}^{-s}C^{4-s}\left\{1+\mcal{O}\left(\frac{1}{\ln T}\right)\right\}\\ 
& = \frac{\zeta^4(2\sigma)}{2\pi^2}\sum_{s=0}^4c_sC^{4-s}\left\{1+\mcal{O}_s\left(\frac{1}{\ln T}\right)\right\}, 
\end{split}
\ee 
where 
\be \label{3.4} 
c_s=b_s\{\zeta(2\sigma)\}^{-s},\ s=0,1,\dots,4,\ c_0=1, 
\ee  
and 
\be \label{3.5} 
\left\{1+\mcal{O}\left(\frac{1}{\ln T}\right)\right\}^{4-s}=1+\mcal{O}_s\left(\frac{1}{\ln T}\right),\ \left|\mcal{O}_s\left(\frac{1}{\ln T}\right)\right|<\frac{A_s(\sigma)}{\ln T}. 
\ee  

\subsection{} 

Now, we put 
\be \label{3.6} 
\frac{\sum_{s=0}^4c_sC^{4-s}\left\{1+\mcal{O}_s\left(\frac{1}{\ln T}\right)\right\}}{\sum_{s=0}^4 c_sC^{4-s}}=1+\frac{\sum_{s=0}^4 c_sC^{4-s}\mcal{O}_s\left(\frac{1}{\ln T}\right)}{\sum_{s=0}^4 c_sC^{4-s}}. 
\ee  
Further, we also use the formula\footnote{See (\ref{3.2}).}
\be \label{3.7} 
C\sim \frac{\ln T}{\zeta(2\sigma)},\ T\to\infty. 
\ee 
It is then true that\footnote{See (\ref{3.4}) -- (\ref{3.6}), (\ref{3.7}).} 
\be \label{3.8} 
\begin{split}
& \frac{\sum_{s=0}^4 c_sC^{4-s}\mcal{O}_s\left(\frac{1}{\ln T}\right)}{\sum_{s=0}^4 c_sC^{4-s}}= 
\frac{\sum_{s=0}^4 c_sC^{-s}\mcal{O}_s\left(\frac{1}{\ln T}\right)}{\sum_{s=0}^4 c_sC^{-s}}= \\ 
& \frac{\mcal{O}_0\left(\frac{1}{\ln T}\right)+\sum_{s=1}^4 c_sC^{4-s}\mcal{O}_s\left(\frac{1}{\ln T}\right)}{1+\sum_{s=1}^4 c_sC^{-s}}< \\ 
& \frac{1}{\ln T}\frac{B_0+\frac{B_1}{\ln T}+\dots+\frac{B_4}{\ln T}}{1-\frac{D_1}{\ln T}-\dots-\frac{D_4}{\ln T}}<\frac{1}{\ln T}\frac{B_0+\frac 12}{1-\frac{1}{2}}=\frac{2B_0+1}{\ln T}, 
\end{split}
\ee 
where 
\bdis 
B_s,D_s>0,\ s=1,\dots,4,\ B_0=A_0. 
\edis 
Now, it follows from (\ref{3.6}) by (\ref{3.8}), that 
\be \label{3.9} 
\begin{split}
& \sum_{s=0}^4c_sC^{4-s}\left\{1+\mcal{O}_s\left(\frac{1}{\ln T}\right)\right\}= \\ 
& \left\{1+\mcal{O}\left(\frac{1}{\ln T}\right)\right\}\sum_{s=0}^4 c_sC^{4-s}, 
\end{split}
\ee 
and subsequently we obtain the next continuation of (\ref{3.3}) 
\be \label{3.10} 
= \frac{\zeta^4(2\sigma)}{2\pi^2}\left\{1+\mcal{O}\left(\frac{1}{\ln T}\right)\right\}\sum_{s=0}^4 c_sC^{4-s}. 
\ee 

\subsection{} 

It is also true, that\footnote{See (\ref{1.1}).} 
\be \label{3.11} 
\begin{split}
& \int_0^T\left|\zf\right|^4{\rm d}t+\mcal{O}(T^{7/8+\epsilon})= \\ 
& \left\{1+\mcal{O}\left(\frac{T^{-1/8+\epsilon}}{\ln^4 T}\right)\right\}\int_0^T\left|\zf\right|^4{\rm d}t. 
\end{split}
\ee 
As a consequence, we have the following result.\footnote{See (\ref{1.1}), (\ref{3.1}), (\ref{3.3}), (\ref{3.10}) and (\ref{3.11}).} 

\begin{mydef51}
\be \label{3.12} 
\begin{split}
& \int_0^T\left|\zf\right|^4{\rm d}t= \\ 
& \frac{\zeta^4(2\sigma)}{2\pi^2}T\left\{1+\mcal{O}\left(\frac{1}{\ln T}\right)\right\}\sum_{s=0}^4 c_s
\left\{
\frac{\int_T^{\overset{1}{T}}|\zf|^2{\rm d}t}{\int_T^{\overset{1}{T}}|\zeta(\sigma+it)|^2{\rm d}t}
\right\}^{4-s}, 
\end{split}
\ee 
where 
\be \label{3.13} 
c_0=1,\ c_s=2\pi^2\{\zeta(2\sigma)\}^{-s}a_s,\ s=1,2,3,4, 
\ee  
for every fixed $\sigma\geq \frac 12+\epsilon$. 
\end{mydef51}

\subsection{} 

In the next parts of this work we use the formula (\ref{3.12}) in the following form. 
\begin{mydefII1}
\be \label{3.14} 
\begin{split}
& \left\{\int_T^{\overset{1}{T}}|\zeta(\sigma+it)|^2{\rm d}t\right\}^4 \times 
\left\{ \int_0^T\left|\zf\right|^4{\rm d}t \right\} \times \\ 
& \left(
\sum_{s=0}^4 c_s\left\{\int_T^{\overset{1}{T}}|\zf|^2{\rm d}t\right\}^{4-s}\times 
\left\{\int_T^{\overset{1}{T}}|\zeta(\sigma+it)|^2{\rm d}t\right\}^s
\right)^{-1}= \\ 
& \frac{\zeta^4(2\sigma)}{2\pi^2}T\left\{1+\mcal{O}\left(\frac{1}{\ln T}\right)\right\}, \ 
\sigma\geq \frac 12+\epsilon,\ T\to\infty. 
\end{split}
\ee 
\end{mydefII1} 

\section{New $\zeta$-functional and corresponding equivalent of the Fermat-Wiles theorem} 

\subsection{} 

Next, the substitution 
\be \label{4.1} 
T=\frac{2\pi^2}{\zeta^4(2\sigma)}x\tau,\ x>0,\ \{T\to\infty\} \Leftrightarrow \{\tau\to\infty\}
\ee 
into (\ref{3.14}) gives the following functional (as the cross-breed of of the Ingham formula and D. R. Heath-Brown formula). 

\begin{mydef11}
\be \label{4.2} 
\begin{split} 
& \lim_{\tau\to\infty}\frac{1}{\tau} \left\{\int_{\frac{2\pi^2}{\zeta^4(2\sigma)}x\tau}^{[\frac{2\pi^2}{\zeta^4(2\sigma)}x\tau]^1}|\zeta(\sigma+it)|^2{\rm d}t\right\}^4\times \int_0^{\frac{2\pi^2}{\zeta^4(2\sigma)}x\tau}\left|\zf\right|^4{\rm d}t\times \\ 
& \left( \sum_{s=0}^4c_s\left\{\int_{\frac{2\pi^2}{\zeta^4(2\sigma)}x\tau}^{[\frac{2\pi^2}{\zeta^4(2\sigma)}x\tau]^1}\left|\zf\right|^2{\rm d}t\right\}^{4-s}\right. \times \\ 
& \left. \left\{\int_{\frac{2\pi^2}{\zeta^4(2\sigma)}x\tau}^{[\frac{2\pi^2}{\zeta^4(2\sigma)}x\tau]^1}|\zeta(\sigma+it)|^2{\rm d}t\right\}^s\right)^{-1} = x 
\end{split} 
\ee 
for every fixed $x>0$ and every fixed $\sigma\geq\frac 12+\epsilon$, where 
\bdis 
c_0=1,\ c_s=2\pi^2\{\zeta(2\sigma)\}^{-s}a_s,\ s=1,2,3,4. 
\edis 
\end{mydef11}

In the special case 
\bdis 
x\to\FR,\ x,y,z,n\in\mbb{N},\ n\geq 3 
\edis 
we obtain from (\ref{4.2}) the following lemma. 

\begin{mydef52}
\be \label{4.3} 
= (4.2) \ \mbox{with}\ x=\FR. 
\ee 
\end{mydef52} 

And consequently, we obtain from (\ref{4.3}) the next result. 

\begin{mydef12}
The $\zeta$-condition 
\be \label{4.4} 
\begin{split} 
& \lim_{\tau\to\infty}\frac{1}{\tau} \left\{\int_{\frac{2\pi^2}{\zeta^4(2\sigma)}\FR\tau}^{[\frac{2\pi^2}{\zeta^4(2\sigma)}\FR\tau]^1}|\zeta(\sigma+it)|^2{\rm d}t\right\}^4\times \\ 
& \int_0^{\frac{2\pi^2}{\zeta^4(2\sigma)}\FR\tau}\left|\zf\right|^4{\rm d}t\times \\ 
& \left( \sum_{s=0}^4c_s\left\{\int_{\frac{2\pi^2}{\zeta^4(2\sigma)}\FR\tau}^{[\frac{2\pi^2}{\zeta^4(2\sigma)}\FR\tau]^1}\left|\zf\right|^2{\rm d}t\right\}^{4-s}\right. \times \\ 
& \left. \left\{\int_{\frac{2\pi^2}{\zeta^4(2\sigma)}\FR\tau}^{[\frac{2\pi^2}{\zeta^4(2\sigma)}\FR\tau]^1}|\zeta(\sigma+it)|^2{\rm d}t\right\}^s\right)^{-1} \not= 1 
\end{split} 
\ee 
on the class of all Fermat's rationals represents the next $\zeta$-equivalent of the Fermat-Wiles theorem for every fixed $\sigma>\frac 12+\epsilon$. 
\end{mydef12} 

\section{On complicated structure of our almost linear increments of the Hardy-Littlewood integral (1918)} 

\subsection{} 

Let us remind the first $zeta$-functional (see \cite{9}) 
\be \label{5.1} 
\lim_{\tau\to\infty}\frac{1}{\tau}\int_{\frac{x}{1-c}\tau}^{[\frac{x}{1-c}\tau]^1}\left|\zf\right|^2{\rm d}t=x,\ x>0. 
\ee  

Current result gives the quotient of the two limits (\ref{4.4}) and (\ref{5.1}), i.e. the quotient of the corresponding integrals\footnote{Factors $\frac{1}{\tau}$ are, of course, cancelled.}. 

\begin{mydef13}
The following equivalence holds true: 
\be \label{5.2} 
\begin{split}
& \int_{\frac{x}{1-c}\tau}^{[\frac{x}{1-c}\tau]^1}\left|\zf\right|^2{\rm d}t\sim \\ 
&  \left\{\int_{\frac{2\pi^2}{\zeta^4(2\sigma)}x\tau}^{[\frac{2\pi^2}{\zeta^4(2\sigma)}x\tau]^1}|\zeta(\sigma+it)|^2{\rm d}t\right\}^4\times \int_0^{\frac{2\pi^2}{\zeta^4(2\sigma)}x\tau}\left|\zf\right|^4{\rm d}t\times \\ 
& \left( \sum_{s=0}^4c_s\left\{\int_{\frac{2\pi^2}{\zeta^4(2\sigma)}x\tau}^{[\frac{2\pi^2}{\zeta^4(2\sigma)}x\tau]^1}\left|\zf\right|^2{\rm d}t\right\}^{4-s}\right. \times \\ 
& \left. \left\{\int_{\frac{2\pi^2}{\zeta^4(2\sigma)}x\tau}^{[\frac{2\pi^2}{\zeta^4(2\sigma)}x\tau]^1}|\zeta(\sigma+it)|^2{\rm d}t\right\}^s\right)^{-1} ,\ \tau\to\infty, 
\end{split}
\ee 
for every fixed $x>0$ and every fixed $\sigma\geq\frac{1}{2}+\epsilon$. 
\end{mydef13} 

\begin{remark}
The equivalence (\ref{5.2}) can be viewed as follows: it expresses the asymptotic \emph{decomposition} of our almost linear increments of the Hardy-Littlewood integral (1918) by means of a complicated composition (non-linear and non-local) of the integrals contained in the cross-breed (\ref{4.2}) of the Ingham and D.R. Heath-Brown formula (1926), (1979). 
\end{remark} 

\section{The construction of finite chains of equivalences for our almost linear increments of the Hardy-Littlewood integral (1918)}

\subsection{} 

Let us remind the next set (for example) of our functionals: (see \cite{10}, (1.3), (5.7)) 
\be \label{6.1} 
\lim_{\tau\to\infty}\frac{1}{\tau}\left\{\sum_{\frac{x}{1-c}\tau<n\leq [\frac{x}{1-c}\tau]^1}d(n)\right\}=x, 
\ee 
(see \cite{11}, (3.10), (5.7)) 
\be \label{6.2} 
\lim_{\tau\to\infty}\frac{1}{\tau}\left\{\sum_{\frac{x}{1-c}\tau< t_\nu\leq [\frac{x}{1-c}\tau]^1}\zfn\right\}=\frac{1}{\pi}x, 
\ee 
(see \cite{12}, (4/6)) 
\be \label{6.3} 
\lim_{\tau\to\infty}\ln\left\{\frac{\Gamma([\frac{x}{1-c}\tau]^1)}{\Gamma(\frac{x}{1-c}\tau)}\right\}^{\frac{1}{\tau}} = x, 
\ee  
(see \cite{15}, (3.3)) 
\be \label{6.4} 
\lim_{\tau\to\infty}\frac{1}{\tau}\int_1^{\frac{x}{\zeta(2\sigma)}\tau}|\zeta(\sigma+it)|^2{\rm d}t=x, 
\ee  
and (see \cite{15}, (3.14)) 
\be \label{6.5} 
\lim_{\tau\to\infty}\frac{1}{\tau}\int_0^{\frac{x}{\bar{c}(l)}\tau}|S_1(t)|^{2l}{\rm d}t=, 
\ee  
where 
\be \label{6.6} 
S_1(t)=\frac{1}{\pi}\int_0^t\arg\left\{\zeta\left(\frac 12+iu\right)\right\}{\rm d}u, 
\ee  
for every fixed $x>0$, $l\in\mbb{N}$ and $\sigma\geq\frac{1}{2}+\epsilon$. 

\subsection{} 

Now, if we use consecutively simple operations from the beginning of the section 5.1 on the continuation of (\ref{5.2}) by means of the functionals (\ref{6.1}) -- (\ref{6.5}), then we obtain the following result. 

\begin{mydef14}
It is true that for every fixed 
\bdis 
x>0, \ l\in\mbb{N},\ \sigma\geq\frac 12+\epsilon
\edis 
there is the following finite chain of equivalences 
\be \label{6.7} 
\begin{split}
& \int_{\frac{x}{1-c}\tau}^{[\frac{x}{1-c}\tau]^1}\left|\zf\right|^2{\rm d}t\sim \\ 
& \left\{\int_{\frac{2\pi^2}{\zeta^4(2\sigma)}x\tau}^{[\frac{2\pi^2}{\zeta^4(2\sigma)}x\tau]^1}|\zeta(\sigma+it)|^2{\rm d}t\right\}^4\times \int_0^{\frac{2\pi^2}{\zeta^4(2\sigma)}x\tau}\left|\zf\right|^4{\rm d}t\times \\ 
& \left( \sum_{s=0}^4c_s\left\{\int_{\frac{2\pi^2}{\zeta^4(2\sigma)}x\tau}^{[\frac{2\pi^2}{\zeta^4(2\sigma)}x\tau]^1}\left|\zf\right|^2{\rm d}t\right\}^{4-s}\right. \times \\ 
& \left. \left\{\int_{\frac{2\pi^2}{\zeta^4(2\sigma)}x\tau}^{[\frac{2\pi^2}{\zeta^4(2\sigma)}x\tau]^1}|\zeta(\sigma+it)|^2{\rm d}t\right\}^s\right)^{-1} \sim \\ 
& \sum_{\frac{x}{1-c}\tau<n\leq [\frac{x}{1-c}\tau]^1}d(n)\sim 
\sum_{\frac{x}{1-c}\tau< t_\nu\leq [\frac{x}{1-c}\tau]^1}\pi\zfn\sim 
\frac{\Gamma([\frac{x}{1-c}\tau]^1)}{\Gamma(\frac{x}{1-c}\tau)}\sim \\ 
& \int_1^{\frac{x}{\zeta(2\sigma)}\tau}|\zeta(\sigma+it)|^2{\rm d}t \sim 
\int_0^{\frac{x}{\bar{c}(l)}\tau}|S_1(t)|^{2l}{\rm d}t. 
\end{split}
\ee 
\end{mydef14} 

\begin{remark}
Let us remind that the integral 
\be \label{6.8} 
\int_{\frac{x}{1-c}\tau}^{[\frac{x}{1-c}\tau]^1}\left|\zf\right|^2{\rm d}t
\ee 
plays the crucial role in our theory. Namely, we can consider it as \emph{a basic state} in the chain (\ref{6.7}) of equivalent states. Consequently, if we use the terminology from the quantum mechanics, every other object in the chain (\ref{6.7}) represents degenerate state of the basic state (\ref{6.8}). 
\end{remark}

\section{Existence of a continuum set of different equivalency chains} 

\subsection{} 

Let us remind that the substitution 
\be \label{7.1} 
T=\frac{x}{1-c}\tau,\ x>0,\ T\in (T_0,\infty), 
\ee 
where $T_0>0$ is sufficiently big, corresponds to the set of the functionals (\ref{5.1}), (\ref{6.1}) -- (\ref{6.3}) while the substitutions 
\be \label{7.2} 
T=\frac{x}{\zeta(2\sigma)}\tau,\ T=\frac{x}{\bar{c}(l)}\tau
\ee 
correspond to the functionals (\ref{6.4}) and (\ref{6.5}), respectively. Next, let 
\be \label{7.3} 
x\in [1,\alpha] 
\ee   
for a fixed $\alpha>1$. Then, it follows from (\ref{7.1}), (\ref{7.2}) that 
\bdis 
\begin{split}
& \frac{1-c}{x}T_0\leq (1-c)T_0=\tau_0^1, \\ 
& \frac{\zeta(2\sigma)}{x}T_0\leq \zeta(2\sigma)T_0=\tau_0^2(\sigma), \\ 
& \frac{\bar{c}(l)}{x}T_0\leq \bar{c}(l)T_0=\tau_0^3(l),  
\end{split} 
\edis 
for $x\in [1,\alpha]$.  

\begin{remark}
Of course, the values 
\be \label{7.4} 
\tau\geq\tau_0^4(\sigma,l)=\max\{\tau_0^1,\tau_0^2(\sigma),\tau_0^3(l)\} 
\ee  
are appropriate for every $x\in [1,\alpha]$ and every fixed 
\bdis 
\sigma\geq \frac{1}{2}+\epsilon,\ l\in\mbb{N}. 
\edis 
\end{remark} 

\subsection{} 

Let 
\be \label{7.5} 
x_1,x_2\in [1,\alpha],\ x_1\not=x_2. 
\ee  
Then, since\footnote{See (\ref{5.1}).} 
\be \label{7.6} 
\begin{split}
& \lim_{\tau\to\infty}\frac{1}{\tau}\int_{\frac{x_1}{1-c}\tau}^{[\frac{x_1}{1-c}\tau]^1}\left|\zf\right|^2{\rm d}t=x_1, \\ 
& \lim_{\tau\to\infty}\frac{1}{\tau}\int_{\frac{x_2}{1-c}\tau}^{[\frac{x_2}{1-c}\tau]^1}\left|\zf\right|^2{\rm d}t=x_2, 
\end{split}
\ee   
we have 
\be \label{7.7} 
\lim_{\tau\to\infty}\frac{\int_{\frac{x_2}{1-c}\tau}^{[\frac{x_2}{1-c}\tau]^1}\left|\zf\right|^2{\rm d}t}{\int_{\frac{x_1}{1-c}\tau}^{[\frac{x_1}{1-c}\tau]^1}\left|\zf\right|^2{\rm d}t}=\frac{x_2}{x_1}\not=1. 
\ee  

\begin{remark}
The following equivalence holds true 
\be \label{7.8} 
\begin{split} 
& \{x_1\not=x_2\} \Leftrightarrow \\ 
&  
\left\{  
\int_{\frac{x_1}{1-c}\tau}^{[\frac{x_1}{1-c}\tau]^1}\left|\zf\right|^2{\rm d}t\nsim 
\int_{\frac{x_2}{1-c}\tau}^{[\frac{x_2}{1-c}\tau]^1}\left|\zf\right|^2{\rm d}t,\ \tau\to\infty
\right\}.
\end{split} 
\ee 
\end{remark} 

\subsection{} 

Let the symbol 
\be \label{7.9} 
\Phi(\tau;x,\sigma,l)
\ee  
denote the chain of equivalences (\ref{6.7}). Then, of course, the equality 
\be \label{7.10} 
\Phi(\tau;x,\sigma,l)=\Phi(\tau;x_1,\sigma,l),\ x_1,x_2\in[1,\alpha] 
\ee 
holds true if and only if the first chain in (\ref{7.10}) contains an element that is equivalent with some elements of the second chain. 

\begin{remark}
Now, it is true\footnote{See (\ref{7.8}).}: 
\be \label{7.11} 
x_1\not=x_2 \ \Rightarrow \ \Phi(\tau;x_1,\sigma,l)\not=\Phi(\tau;x_2,\sigma,l). 
\ee 
\end{remark} 

\subsection{} 

Let us remind that the functional (\ref{5.1}) is strictly increasing on the segment $[1,\alpha]$. Consequently, the next corollary follows from (\ref{6.7}) and (\ref{7.11}). 

\begin{mydef41}
The set 
\bdis 
\{\Phi(\tau;x,\sigma,l)\}_{x\in[1,\alpha]}
\edis 
is the continuum set of different chains for every fixed $\sigma\geq \frac{1}{2}+\epsilon$, $l\in\mbb{N}$ and $\alpha>1$.  
\end{mydef41}

I would like to thank Michal Demetrian for his moral support of my study of Jacob's ladders.

\end{document}